\documentclass[12pt,a4paper]{amsart}

\usepackage{amsfonts}
\usepackage{amssymb}
\usepackage{amsmath}   
\usepackage{amsthm}    
\usepackage{amscd}     

\usepackage[english]{babel}

\newcommand{\NN}{{\mathbb N}}
\newcommand{\ZZ}{{\mathbb Z}}

\newcommand{\CC}{{\mathbb C}}

\newcommand{\PP}{{\mathbb P}}    
\newcommand{\GG}{{\mathbb G}}    

\newcommand{\PPn}{{\PP^n}}

\newcommand{\slc}{{\mathcal {SL}}}

\DeclareMathOperator{\HH}{H}
\DeclareMathOperator{\Hom}{Hom}
\DeclareMathOperator{\End}{End}
\DeclareMathOperator{\SL}{SL}
\DeclareMathOperator{\GL}{GL}
\DeclareMathOperator{\Sym}{Sym}
\DeclareMathOperator{\coker}{Coker}
\DeclareMathOperator{\kur}{Kur}
\DeclareMathOperator{\Or}{O}
\DeclareMathOperator{\Id}{Id}

\newcommand{\PPO}{\PP \Or}
\newcommand{\OPN}{{\mathcal O}_{\PPn}}
\newcommand{\FAG}{F_{\alpha,\gamma}}
\newcommand{\QAG}{Q_{\alpha,\gamma}}
\newcommand{\FT}{\widetilde F_o}
\newcommand{\QT}{\widetilde Q_o}
\newcommand{\quot}{{\it Quot}_{{\mathcal V}|\PP^n}}

\newtheorem{thm}{Theorem}[section]
\newtheorem{lemma}[thm]{Lemma}
\newtheorem{prop}[thm]{Proposition}

\newenvironment{remark}{{\bf Remark.}}

\title{Weighted Tango Bundles On $\PPn$ And Their Moduli Spaces}
\author{Paolo Cascini}
\address{Dipartimento di Matematica\\ Viale Morgagni 67 A\\ 50134 Firenze\\ Italy}
\email{cascini@math.unifi.it}

\begin{document}
\maketitle
\begin{abstract} We define a new class of algebraic $(n-1)$-bundles on $\PPn$, that contains the bundles introduced by Tango \cite{Tan} and their stable generalized pull-backs; we show that these bundles are invariant under small deformations and that they correspond to smooth points of moduli spaces.
\end{abstract}

\vskip 1 cm

It is a very difficult problem to find examples of non-splitting algebraic vector bundles on the complex projective space $\PPn$ whose rank is less than $n$. In particular for $n\ge 6$ the only known examples are essentially the mathematical instantons \cite{AO2} (for odd $n$) and the bundles introduced by Tango \cite{Tan}: all of them have rank $n-1$. Of course, pulling back the Tango bundles by a finite morphism $\PPn\longrightarrow \PPn$ gives other examples of rank $n-1$ bundles.

In \cite{Hor}, Horrocks introduced a new technique of constructing new bundles from old ones, which generalizes the pull-back. This method, that we can call {\it generalized pull-back}, has been extensively studied in \cite{AO} and \cite{AO1} and it applies only to bundles whose symmetry group contains a copy of $\CC^*$.

In this paper we show that, for any $n\ge 3$, there exists a Tango bundle that is $\SL(2)$-invariant: hence the generalized pull-back allows us to define a new class of $(n-1)$-bundles on $\PPn$.

More precisely, let $\alpha, \gamma$ be integer numbers  such that $\gamma > n\alpha\ge 0$ and let $\QAG$ be the bundles on $\PPn$ described by the exact sequence:
$$0\rightarrow \OPN(-\gamma)\rightarrow \bigoplus_{k=0}^{n}\OPN((n-2k)\alpha)\rightarrow\QAG\rightarrow 0.$$
$\QAG$ can also be defined as the generalized pull-back of the quotient bundle on $\PPn$ and,
in particular, $Q_{0,1}$ is the quotient bundle.
Let us define the rank $2n-1$ vector bundle:
$${\mathcal V}=S^{2(n-1)}(\OPN(\alpha)\oplus\OPN(-\alpha))=\bigoplus_{k=0}^{2(n-1)}\OPN((2n-1-2k)\alpha).$$ 
It will be proven that there exists an exact sequence of algebraic vector bundles over $\PPn$:
\begin{equation}\label{paolo}
0\rightarrow \QAG(-\gamma)\rightarrow {\mathcal V}\rightarrow \FAG(\gamma)\rightarrow 0.
\end{equation}
The $(n-1)$-bundles $\FAG$ are called {\it weighted Tango bundles of weights $ \alpha$ and $\gamma$} and they are stable if and only if $\gamma>2(n-1)\alpha$.
The bundles $F_{0,1}$ are the classical Tango bundles, moreover the
generalized pull-backs of the Tango bundles are contained in the
sequence (\ref{paolo}).  
The main result of this paper is the following:
\begin{thm}\label{princess}
Let $\FAG^o$ be a stable weighted Tango bundle on $\PPn$ of weights
$\alpha$ and $\gamma$ and let $c_i$ be the $i$-th Chern class of
$\FAG^o$ (in particular $c_1=0$).  
There exists a smooth neighborhood of the point of the moduli space 
${\mathcal M}_{\PPn}(0,c_2,\dots,c_{n-1})$ corresponding to $\FAG^o$
entirely consisting of weighted Tango bundles of weights $\alpha$ and
$\gamma$. 
\end{thm}

I would like to express my gratitude to professor V.~Ancona for his
invaluable guidance and to professor G.~Ottaviani for his insightful
suggestions. 
 
\vskip 1 cm

\section{\bf Introduction.}

Let $V$ be a $(n+1)$-dimensional vector space over  $\CC$, and let $\PPn=\PP(V)$:
it is possible to show (cf. \cite{JSW}) that a Tango bundle $F$ on $\PPn$ is contained in the following exact sequence: 
$$0\rightarrow Q(-1) \rightarrow \frac{\wedge^2V} W\otimes \OPN \rightarrow F(1)\rightarrow 0;$$
here $Q$ is the quotient bundle (cf. \cite{OSS}) on $\PPn$ and $W\subseteq \wedge^2 V$ is a linear subspace such that:
\begin{equation}\label{ciao}
\begin{cases}
\dim_{\CC}\PP(W)  =m-1 \\ 
\PP(W)\cap \GG(1,n) =\emptyset
\end{cases}
\end{equation}
where $m= {\frac{(n-2)~(n-1)} 2}$ and
$\GG(1,n)$ is the Grassmannian of the lines in $\PPn=\PP(V)$: hence $W$ does not contain any decomposable bivectors.
Moore \cite{Moo}  has shown that $F$ is uniquely determined by the subspace $W\subseteq \wedge^2V$ and so by a point of the variety $\GG(m-1,N-1)$, with $N= {\frac{n~(n+1)}2}$; furthermore if $W$ is invariant under the action of a group $G\subseteq\PP\GL(n+1)$ then the Tango bundle, associated to $W$, is $G$-invariant too, i.e. $G\subseteq \Sym F$.  

\vskip 1 cm
\section{\bf Action of $\SL(2)$.}

Let $U$ be a $2$-dimensional vector space over $\CC$ and let us consider the complex projective space $\PPn=\PP(S^nU)$: in this way, we have a natural action of  $\SL(2) = \SL(U)$ over $\PPn$.

We want to find a subspace $W\subseteq\wedge^2 S^nU$, $\SL(2)$-invariant and that satisfies (\ref{ciao}). For this purpose we prove the following:

\begin{prop}\label{daphne} The decomposition of $\wedge^2S^nU$ into irreducible representations is given by $S^{2(n-1)}U\oplus S^{2(n-3)}U\oplus S^{2(n-5)}U\oplus\dots$;
moreover if
$W=S^{2(n-3)}U\oplus S^{2(n-5)}U\oplus\dots$,
then $W$ satisfies (\ref{ciao}).
\end{prop}

This proposition immediately implies that for any $n\in \NN$, such a subspace $W$ defines a $\SL(2)$-invariant Tango bundle $F$ on $\PPn$, which is
 described by the exact sequence:
$$0\rightarrow Q(-1) \rightarrow S^{2(n-1)}U\otimes \OPN \rightarrow F(1)\rightarrow 0.$$

Before proceeding with the proof of the proposition, we prove the following lemma: 
\begin {lemma}
Let $\{ v_0,\dots,v_n\} $ be a basis of $V$ and $\omega\in\GG (1,n)\subseteq\wedge^2V$ a non-vanishing decomposable bivector, then:
$$\omega = x_{i_0,j_0}(v_{i_0}\wedge v_{j_0})+ \sum_{i+j>i_0+j_0}x_{i,j}(v_i\wedge v_j)$$
where $x_{i,j}\in\CC$ and $x_{i_0,j_0}\neq 0$.
\end {lemma}

\noindent
\begin{remark} In order to simplify the notations, we will often write $v_{i,j}$ instead of $v_i\wedge v_j$.
\end{remark}
\begin{proof}
We proceed by induction on $n$.
For $n=1$, there is nothing to prove.

Let us suppose now $n>1$ and let $\omega =v\wedge v'$ where $v=\sum
x_i v_i$ and $v'=\sum y_i v_i$. 

Let $z_{i,j}=x_iy_j-x_jy_i$ then 
$$\omega=\sum^{i<j}_{i+j\ge k_0} z_{i,j}v_{i,j},$$
where $k_0=\min\{k|z_{i,j}=0 \text{ if } i+j=k\}$.

If there exist $i_0,j_0\neq 0$ such that $i_0+j_0=k$, and
$z_{i_0,j_0}\neq 0$ then, since
$z_{0,i_0}=z_{0,j_0}=0$, it easily follows $x_0=y_0=0$: thus the lemma
is true by induction.

Otherwise, if such $i_0,j_0$ do not exist, then:
$$\omega=z_{0,k}v_{0,k}+\sum^{i<j}_{i+j>k_0} z_{i,j}v_{i,j}.$$ 
\end{proof}

\vskip .4 cm 
\begin{proof}[Proof of proposition \ref{daphne}]

\

Let $V=S^nU$ and let $\{x,y\}$ be a basis of $V$: if $v_0=x^n,\dots, v_n=y^n$, then $\{v_0\dots v_n\}$ is a basis of $V$.
The weights of $S^nU$ are $\{n, n-2,\dots,-n\}$ (cf. \cite{Ful}, pag. 146--153) and since the weights of $\wedge^2 S^nU$ are given by the sums of couples of different weights of $S^nU$, it easily follows:
$$\wedge^2S^nU = S^{2(n-1)}U\oplus S^{2(n-3)}U\oplus S^{2(n-5)}U\oplus\dots$$
Indeed if $W= S^{2(n-3)}U\oplus S^{2(n-5)}U\oplus\dots$, then $\dim_{\CC} W=m$.
 
Let us prove now that $W$ does not contain any decomposable bivector, as required.
We suppose that there exists $\omega \in W\cap \GG(1,n)$, such that $\omega\neq 0$; by the previous lemma, we get:
$$\omega = x_{i_0,j_0}~v_{i_0,j_0}+ \sum_{i+j>i_0+j_0}x_{i,j}~v_{i,j}$$
where $x_{i_0,j_0}\neq 0$. We want to show that, in this case,  there
exists a vector of weight $2(n-1)$ in $W$: this contradicts with the
fact that $S^{2(n-1)}U\cap W = \{0\}$.  

Let 
$Y=\begin{pmatrix}0 & 0\cr 1 & 0\end{pmatrix}, H=\begin{pmatrix}1 & 0\cr 0 &-1\end{pmatrix}\in \slc(2)$,
and let $\tilde Y, \tilde H$ be the corresponding endomorphisms of $\wedge^2 S^nU$.
If we suppose $v_{n+1}=0$, we have: 
$$\tilde Y(v_{i,j})=(n-i)~v_{i+1,j}+(n-j)~v_{i,j+1}\qquad \text{for any }i,j=0,\dots,n.$$
Hence if $k=(2n-1)-i_0-j_0$, then $x_{i_0,j_0}~\tilde Y^{(k)}(v_{i_0}\wedge v_{j_0})=\tilde Y^{(k)}(\omega)\in W$. On the other hand it results that  $ \tilde Y^{(k)}(v_{i_0,j_0})=m~v_{n,n-1}$, where $m$ is a positive integer: this implies 
that $v_{n,n-1}\in W$ and since $\tilde H (v_{n,n-1})=-2(n-1)v_{n,n-1} $, we see that $W$ contains a vector of weight $2(n-1)$.
\end{proof}

\begin{remark}
Moore \cite{Moo} has shown that the Tango bundles on $\PP^4$ have all
symmetry groups isomorphic to $\PPO(3)$ and that $\PP\GL(5)$ acts
transitively on the moduli space of the Tango bundles ${\mathcal
M}_{\PP^4}(0,2,2)$. 
In higher dimensions the situation is different: in fact, with the
help of the software Macaulay 2 \cite{Mac}, 
it has been possible to prove that on $\PP^5$ the generic Tango bundle
has a discrete symmetry group  and that there exist Tango bundles with
the symmetry group isomorphic to $\CC^*$ (for istance the one defined
by
$W=<v_{0,5}+5v_{2,3},v_{1,4}+3v_{2,3},v_{0,4}-2v_{1,3},v_{2,5}+v_{3,4},
v_{0,3}+3v_{1,2},2v_{2,5}-3v_{3,4}>$).  

The algorithm needed to calculate the dimension of the orbit of a
subspace $W_0\subseteq\wedge^2 V$ (where $n=5$) 
under the action of $\PP \GL(6)$ was communicated to the author by
G. Ottaviani. We describe the 
fundamental steps of it: 

\begin{enumerate}

\item Let us choose $m$ as a $(6\times 15)$-matrix whose rows represent the generators of
the  subspace $W_0\subseteq\wedge^2 V$;

\item We denote by $g = \{g_{i,j}\}$ a generic $(6\times
6)$-matrix and let's  define $m'=m*\wedge^2 g$:

\end{enumerate}

\noindent $m'$ represents the image $g W_0$  of the matrix $g\in  \PP
\GL(6)$ by the map  $\eta:\PP \GL(6)\rightarrow \GG(6,\wedge^2 V)$;
By the Plucker embedding $\phi:\GG(m,\wedge^2 V)\hookrightarrow \PP^{5004}$, 
the dimension of the orbit of $W_0$ is equal to the dimension of the
ideal generated by the minors $6\times 6$ of $m'$, but its calculation 
is, computationally, too difficult. 
Therefore in order to make the computation easier, we first calculate the derivative $d(\phi\circ  \eta)$ 
at the identity matrix and then we compute
the dimension of its image: this number is exactly the dimension
of the orbit. We proceed as follows: 

\begin{enumerate}
\stepcounter{enumi}
\stepcounter{enumi}
\item Let $v_1(g),\dots,v_6(g)$ be the rows of $m'$, and let 
$v_i(g)_{g_{i,j}}= \frac {\partial v_i(g)} {\partial g_{i,j}}$. 

\end{enumerate}

\noindent In order to compute the derivative $d(\phi\circ  \eta)$,
we remind that, for any $I\subseteq\{1,\dots,15\}$ such that $\#I=6$,
we have:
$$\frac \partial {\partial g_{i,j}} \det \begin{pmatrix}v_0^I(g)\cr
\vdots\cr v_6^I(g)\end{pmatrix} = 
\det \begin{pmatrix}v_0^I(g)_{g_{i,j}}\cr 
v_1^I(g)\cr \vdots\cr v_6^I(g)\end{pmatrix}
+ \dots +
\det \begin{pmatrix}v_0^I(g)\cr\vdots\cr 
v_5^I(g) \cr v_6^I(g)_{g_{i,j}}\end{pmatrix}
$$   
where $v_i^I(g)$ denotes the vector composed by the components of
$v_i(g)$ with index in $I$.

\begin{enumerate}
\stepcounter{enumi}
\stepcounter{enumi}
\stepcounter{enumi}
\item Let's define $M_{i,j}^k = \begin{pmatrix}v_1(\Id_6)\cr\vdots\cr
v_k(\Id_6)_{g_{i,j}} \cr\vdots\cr v_6(\Id_6)\end{pmatrix}$;

\item let
$p_{i,j}$ be the sum of the vectors in $\PP^{5004}$ defined by the minors
of $M_{i,j}^k$ with $k=1,\dots,6$;

\item The rank of the matrix $\begin{pmatrix}p_{1,1}\cr
p_{1,2}\cr\vdots \cr p_{6,6}\end{pmatrix}$ is the dimension of the
orbit of $W_0$. 
\end{enumerate}

\end{remark}

\vskip 1 cm
\section{\bf Weighted Tango Bundles.}

We have shown that for any $n$, there exists a Tango bundle $F$ on $\PP(S^nU)$ that is invariant under the $\CC^*$-action defined by:
$$\begin{pmatrix} t^n &         &       &\cr
                      & t^{n-2} &       &\cr
                      &         &\ddots &\cr
                      &         &       & t^{-n} \end{pmatrix}\in\PP\GL(n+1)
\qquad\text{for any }t\in\CC^*  $$
      
This map induces an embedding of $\CC^*$ in $\Sym F$ and so it is possible to study the pull-backs over $\CC^{n+1}\setminus 0$ of such bundles (cf. \cite{AO, AO1}).

Let us fix $\alpha, \gamma\in \NN$ such that $\gamma>n\alpha$ and let 
$f_0,\dots,f_n\in\CC[x_0,\dots,x_n]$ homogeneous polynomial of degree: 
$$\deg f_k = \gamma + (n -2k)~\alpha \qquad\text{for each }k=0,\dots,n $$
and without common roots.

Let $\phi=(f_0,\dots,f_n)$ and let us take into account the following diagram:
$$\begin{CD}
    \CC^{n+1}\setminus 0   @>\phi>>    S^nU\setminus 0\\
        @V\pi_1 VV                           @VV\pi_2 V \\
       \PP^n            &&             \PP^n
\end{CD}$$

According to \cite{AO, Hor}, there exists an algebraic vector bundle $\FAG$ on $\PPn$ such that $\pi_1^*\FAG=\phi^*\pi_2^* F$. Furthermore, since $Q$ is an homogeneous bundle \cite{OSS}, there exists $\QAG$ such that $\pi_1^*\QAG=\phi^*\pi_2^* Q$.  
Such a bundle is contained in the weighted Euler sequence:
\begin{equation}\label{eulpes}
0\rightarrow \OPN(-\gamma)\rightarrow S^n{\mathcal U}\rightarrow\QAG\rightarrow 0
\end{equation}
where ${\mathcal U}=\OPN(-\alpha)\oplus\OPN(\alpha)$.
In general, we will call {\it weighted quotient bundle of weights $\alpha$ and $\gamma$} any bundles $\QAG$ contained in a sequence (\ref{eulpes}).

On the other hand $\FAG$ is contained in the exact sequence:
\begin{equation}\label{tanpes}
0\rightarrow \QAG(-\gamma)\rightarrow {\mathcal V}\rightarrow \FAG(\gamma)\rightarrow 0
\end{equation}
where ${\mathcal V}=S^{2(n-1)}{\mathcal U}$ and $\QAG$ is the  pull-back over $\CC^{n+1}\setminus 0$ of the quotient bundle $Q$ defined by the map $\phi$.
Also in this case, we will call {\it weighted Tango bundle of weights $\alpha$ and $\gamma$} any bundles $\FAG$ contained in the sequence (\ref{tanpes}), where $\QAG$ is any  weighted quotient bundle of weights $\alpha$ and $\gamma$.

By these sequences, it immediately follows that $c_1(\FAG)=0$ and that $c_i(\FAG)=c_i(\alpha,\gamma)$ for any $i=2,\dots,n-1$ (i.e. the Chern classes do not depend on the map $\phi$).

\begin{prop} A weighted Tango bundle $\FAG$ is stable if and only if  $\gamma>2(n-1)\alpha$.
\end{prop}

\begin{proof} Let $\gamma>2(n-1)\alpha$. By the Hoppe criterion \cite{Ho}, it suffices to show that $\HH^0(\wedge^q\FAG)=0$ for any $q=1,\dots,n-2$.
By the sequence:
$$0\rightarrow S^{k-1}S^n{\mathcal U}(-\gamma)\rightarrow S^{k}S^n{\mathcal U}\rightarrow S^k\QAG\rightarrow 0$$
obtained raising the sequence (\ref{eulpes}) to the $k$-th symmetric power, we see that:
$\HH^i(S^k\QAG(t))=0$ for any $i=1,\dots,n-2$ and $t\in\ZZ$.

On the other hand by (\ref{tanpes}), we have the long exact sequence:
$$0\rightarrow S^q\QAG(-q\gamma)\rightarrow \dots \rightarrow S^k\QAG(-k\gamma)\otimes \wedge^{q-k} {\mathcal V}\rightarrow\dots$$
$$\dots\rightarrow\QAG(-\gamma)\otimes\wedge^{q-1}{\mathcal V}\rightarrow\wedge^q{\mathcal V}\rightarrow\wedge^q\FAG(q\gamma)\rightarrow 0$$
This sequence immediately implies that $\HH^0(\wedge^q\FAG)\subseteq\HH^0(\wedge^q {\mathcal V}(-q\gamma))$, and since
$$\max\{t\in\ZZ|\OPN(t)\subseteq \wedge^q {\mathcal V}(-q\gamma)\}=q((2n-q-1)\alpha-\gamma)<0$$
we have that $\HH^0(\wedge^q\FAG)=0$ for any $q=1,\dots,n-2$, and so $\FAG$ is stable.

Let us prove now that the condition is necessary.
By the sequences:
$$0\rightarrow \OPN(-3\gamma)\rightarrow S^n{\mathcal U}(-2\gamma)\rightarrow \QAG(-2\gamma)\rightarrow 0$$
$$0\rightarrow \QAG(-2\gamma)\rightarrow {\mathcal V}(-\gamma)\rightarrow \FAG\rightarrow 0$$
it follows that if  $\gamma\le 2(n-1)\alpha$, then $\HH^0(\FAG)\neq 0$ and so $\FAG$ cannot be stable.  
\end{proof}

\vskip 1 cm
\section{\bf Small deformations of $\FAG$.}
Let $E$ be a vector bundle on $\PPn$: we will indicate with $(\kur E,e)$ the Kuranishi space of $E$ (cf. \cite {FK}), where $e\in \kur E$ is the point corresponding to the bundle $E$.

We are finally ready to introduce the main result of this paper:

\begin{prop}\label{teo}
Let $\FAG^o$ be a weighted Tango bundle of weights $\alpha$ and $\gamma$.
Every small deformation of $\FAG^o$ is still a weighted Tango bundle and its Kuranishi space is smooth at the point corresponding to $\FAG^o$.
\end{prop}
\noindent

\vskip 0.5 cm
Before proceeding with the proof of the proposition, let us look at some preliminaries:

\begin{lemma}\label{koonteng}
Let $\QAG^o$ be a weighted quotient bundle.
Every small deformation of $\QAG^o$ is still a weighted quotient bundle and the Kuranishi space of $\QAG^o$ is smooth at the point corresponding to its isomorphism class.
\end{lemma}

\begin{proof}
The proof of this lemma is very similar to the proof of prop.~3.1 of \cite{AO}.
\end{proof}

\begin{lemma}\label{l1}
Let $\FAG$ and $\FAG'$ be two isomorphic weighted Tango bundles, defined by the sequences: 
$$0\rightarrow \QAG(-\gamma) \rightarrow {\mathcal V}\rightarrow \FAG(\gamma) \rightarrow 0$$ 
$$0\rightarrow \QAG'(-\gamma) \rightarrow {\mathcal V}\rightarrow \FAG'(\gamma) \rightarrow 0$$ 
where $\QAG$ and $\QAG'$ are weighted quotient bundles.
Then $\QAG$ and $\QAG'$ are isomorphic.
\end{lemma}
\begin{proof}
By joining together the sequences (\ref{eulpes}) and (\ref{tanpes}), we get:
$$0\rightarrow \OPN(-2\gamma){\stackrel{\phi}\rightarrow} S^n{\mathcal U}(-\gamma)\rightarrow {\mathcal V}\rightarrow\FAG(\gamma)\rightarrow 0.$$
By proposition 1.4 of \cite{BS} and by the fact that $-2\gamma<-\gamma-n\alpha$, the last sequence is the minimal resolution of $\FAG(\gamma)$: hence $\QAG(-2\gamma)=\coker \phi$ is directly defined by this resolution. 
\end{proof}

\begin{lemma} \label{l2}
Every isomorphism between two weighted Tango bundles $\FAG \rightarrow \FAG'$ is induced by an isomorphism of sequences:
$$\begin{CD}
0@>>> \QAG(-\gamma) @>>> {\mathcal V} @>>> \FAG(\gamma)  @>>> 0 \\
    &&   @VVV                @VVV               @VVV            \\
0@>>> \QAG'(-\gamma) @>>> {\mathcal V} @>>> \FAG'(\gamma)  @>>> 0 \\
\end{CD}$$
\end{lemma}

\begin{proof}
By the sequence
$$0\rightarrow \OPN(-2\gamma)\otimes {\mathcal V}\rightarrow S^n{\mathcal U}(-\gamma)\otimes {\mathcal V}\rightarrow\QAG(-\gamma)\otimes {\mathcal V}\rightarrow 0,$$
and since 
$$h^1(S^n{\mathcal U}(-\gamma)\otimes {\mathcal V})=h^2(\OPN(-2\gamma)\otimes {\mathcal V})=0,$$
we get $h^1(\QAG(-\gamma)\otimes {\mathcal V})=0;$ hence the lemma is proven.
\end{proof}

\begin{lemma} \label{l3}
Two morphisms $f,f'\in \Hom(\QAG(-\gamma),{\mathcal V})$ give the same element of $\quot$ if and only if there exists an invertible $h\in\End(\QAG(-\gamma))$ such that
$$\ f=f'\circ h.$$
\end{lemma}
 
\begin{proof}
It follows from the definition of $\quot$, (cf. \cite{Hu}).
\end{proof}

\

\begin{proof}[Proof of proposition \ref{teo}]

\

For brevity's sake, we will write $\FT$ instead of $\FAG^{o} $ and $\QT$ for $\QAG^{o}$.
Let also $\sigma_0\in \Hom(\QT(-\gamma),{\mathcal V})$ be such that $\FT=\coker \sigma_0$.

Let ${\mathcal Q}$ be the sub-variety of the irreducible component of $\quot$ composed by all the quotients of the maps $0\rightarrow \QAG(-\gamma){\stackrel {\sigma}\longrightarrow} {\mathcal V}$ for some weighted bundle $\QAG$ and containing the point $\sigma_0$ corresponding to $\FT$: the morphisms $\Phi:({\mathcal Q},\sigma_0)\longrightarrow(\kur {\QT},q_0)$ and $\Psi:({\mathcal Q},\sigma_0)\longrightarrow (\kur {\FT},f_0)$ are canonically defined.

A generic fiber of $\Phi$ is given by all the cokernels of the morphisms  $\QAG(-\gamma)\rightarrow {\mathcal V}$ with a fixed $\QAG$, and so, by lemma \ref{l3}, its dimension is constantly equal ($\alpha$ and $\gamma$ are fixed) to 
$h^0(\QAG^*(\gamma)\otimes {\mathcal V})-h^0(\End \QAG)$.
Hence, since lemma \ref{koonteng} implies that $\dim_{q_0}(\kur {\QT})=h^1(\End \QT)$, we get:
$$\dim_{\sigma_0}{\mathcal Q}=h^0(\QT^*(\gamma)\otimes {\mathcal V})-h^0(\End \QT)+ h^1(\End \QT)$$

Let us study now the morphism $\Psi:\mathcal Q \longrightarrow \kur {\FT}$:
if $\Sigma=\{ \sigma\in \quot| F_{\sigma} \simeq \FT \}$,
then it results $\Psi^{-1}(f_0)\subseteq \Sigma$ and by lemma \ref{l1}, \ref{l2} and \ref{l3}, it follows:
$$\dim_{\sigma_o}\Sigma = h^0(\End {\mathcal V})-\dim \{\varphi \in \End {\mathcal V}|\varphi\cdot\sigma_0 =\sigma_0 \}-h^0(\End \QT).$$
By the sequence:
$$0\rightarrow\FT^* (-\gamma)\otimes {\mathcal V} \rightarrow \End {\mathcal V} \rightarrow \QT^*(\gamma)\otimes {\mathcal V} \rightarrow 0$$
obtained tensoring the dual sequence of (\ref{tanpes}) with ${\mathcal V}$, we have that:
$$\dim \{\varphi \in \End {\mathcal V}|\varphi\cdot\sigma_0 =\sigma_0\}=h^0(\FT^*(-\gamma) \otimes {\mathcal V})$$
and so:
$$\dim_{\sigma_0} \Psi^{-1}(f_0)\le \dim_{\sigma_o}\Sigma = h^0(\End {\mathcal V})-h^0(\FT^* (-\gamma)\otimes {\mathcal V})-h^0(\End \QT).$$
Hence:
$$h^1(\End \FT)\ge \dim_{f_0}(\kur {\FT})\ge h^1(\End \QT)+h^1(\FT^* (-\gamma)\otimes {\mathcal V}).$$
To prove the proposition it suffices to show that $$h^1(\End \FT)\le h^1(\End \QT)+h^1(\FT^* (-\gamma)\otimes {\mathcal V}).$$ 
In fact this implies that $h^1(\End \FT)=\dim_{f_0}(\kur {\FT})$, i.e. $\kur {\FT}$ is smooth at the point $f_0$, and that $\dim_{f_0}(\kur {\FT})= \dim_{\sigma_0}{\mathcal Q}-\dim \Psi^{-1}(f_0)$, i.e. $\Psi$ is surjective.

By the exact sequence:
$$0\rightarrow \QT(-2\gamma)\otimes \FT^*\rightarrow \QT(-\gamma)\otimes {\mathcal V}\rightarrow \End \QT\rightarrow 0$$
and by the vanishing of $\HH^1(\QT(-\gamma)\otimes {\mathcal V})$ and $\HH^2(\QT(-\gamma)\otimes {\mathcal V})$, we have  that $\HH^1(\End \QT)=\HH^2(\QT(-2\gamma)\otimes \FT^*)$.
Hence by the sequence: 
$$0\rightarrow \QT(-2\gamma)\otimes \FT^*\rightarrow \FT^*(-\gamma)\otimes {\mathcal V}\rightarrow \End \FT\rightarrow 0$$
and for what we have seen, we get the sequence of cohomology groups:
$$\dots\rightarrow \HH^1(\FT^*(-\gamma)\otimes {\mathcal V})\rightarrow\HH^1(\End \FT)\rightarrow \HH^1(\End \QT)\rightarrow\dots$$
In particular
$h^1(\End \FT)\le h^1(\End \QT)+h^1(\FT^* (-\gamma)\otimes {\mathcal V})$, as required.
\end{proof}

\vskip 1 cm

Theorem \ref{princess} easily follows from the previous proposition.
In fact if $\gamma\ge 2(n-1)\alpha$, we can consider the canonical algebraic map $\mathcal Q\longrightarrow \mathcal M(0,c_2,\dots,c_{n-1})$. The image of this map is a smooth quasi projective set composed uniquely by weigthed Tango bundles and it is an open neighborhood of $\FAG^o$ in $\mathcal M(0,c_2,\dots,c_{n-1})$.

\vskip 1 cm
\end{document}